\documentclass[10pt,oneside,disablejfam,notitlepage]{article}
\usepackage{framed}
\usepackage[mathscr]{eucal}
\usepackage[dvipdfm]{graphicx}
\usepackage{makeidx}
\usepackage{bm}
\usepackage{amsmath,amssymb,amsthm,amscd}
\usepackage{tabularx}
\usepackage{slashbox}
\usepackage{multicol}
\usepackage{multirow}
\usepackage{hhline}
\usepackage{ccaption}
\usepackage[dvipdfm]{color}

\theoremstyle{definition}
\newtheorem{thm}{Theorem}[section]

\newtheorem{lem}[thm]{Lemma}
\newtheorem{prop}[thm]{Proposition}

\newtheorem{rem}[thm]{Remark}

\begin{document}

\vspace{-1cm}

\begin{center}

\large

\textbf{Non-existence of elliptic curves with everywhere good reduction over some real quadratic fields}

\normalsize

\vspace{3mm}

Shun'ichi Yokoyama and Yu Shimasaki

\vspace{1mm}

(Kyushu University)


\end{center}

\vspace{3mm}

\noindent
\textbf{Abstract} \ \ We prove the non-existence of elliptic curves having good reduction everywhere over some real quadratic fields.

\vspace{4mm}

\section{Introduction}

Throughout this paper, let $K_m$ be the real quadratic field $\mathbb{Q}(\sqrt{m})$ where $m$ is a square-free positive integer with $m \le 100$ and $\mathcal{O}_{K_m}$ the ring of integers of $K_m$. We already know the following results concerning elliptic curves with everywhere good reduction over real quadratic fields (\cite{CL00}, \cite{Ish86}, \cite{Kag98}, \cite{Kag99}, \cite{Kag01}, \cite{Kag11}, \cite{Kid99}, \cite{Kid01}, \cite{KK97}, \cite{Pin82}, \cite{TT11}) :

\begin{thm}
\begin{enumerate}
\item \textit{There are no elliptic curves with everywhere good reduction over $K_m$ if} \\ \textit{$m=2,3,5,10,11,13,15,17,19,21,23,30,31,34,35,39,42,47,53,55,57,58,$}
\\ \textit{$61,66,69,70,73,74,78,82,83,85,89,93,94,95$ and $97$.}
\item \textit{The elliptic curves with everywhere good reduction over $K_m$ are determined completely for $m=6,7,14,22,29,33,37,38,41,65$ and $77$.}
\item \textit{There are elliptic curves with everywhere good reduction over $K_m$ if $m=26,79$ and $86$} (cf.\ \cite{Ish86}, \cite{MSZ89} and Cremona's table \cite{CreTb}).
\end{enumerate}
\end{thm}

In this paper, we prove the non-existence of elliptic curves with everywhere good reduction over three real quadratic fields not appearing in Theorem 1.1. Here is the main theorem:

\begin{thm}
\textit{If $m=43,46$ and $59$, there are no elliptic curves with everywhere good reduction over $K_m$.}
\end{thm}

The following cases are still unknown whether an elliptic curve with everywhere good reduction over $K_m$ exists or not:
\[
m=51, 62, 67, 71, 87, 91.
\]
For three of them, we prove the following conditional result:

\begin{thm}
\textit{If $m=62,67$ and $71$, there are no elliptic curves with everywhere good reduction over $K_m$ which have cubic discriminant.}
\end{thm}

\begin{rem}
In some cases (e.g. $m=77$, cf.\ \cite{Kag11}), we can prove that there exists an elliptic curve with everywhere good reduction over $K_m$ having cubic discriminant.
\end{rem}

Our strategy for the proof is close to that of T. Kagawa \cite{Kag11}. However, we use different kinds of computer softwares and computational techniques. In \cite{Kag11}, all computations were carried out by using KASH \cite{KASH} and SIMATH \cite{SIMATH}. Unfortunately, development of SIMATH had already stopped and some fatal bugs (Tate's algorithm over number fields, for example) remain even now. Thus we switched environment of computation completely and started from a check experiment of Kagawa's results by using Magma \cite{Magma}, Pari-GP \cite{Pari} and Sage \cite{Sage}.

\vspace{4mm}

\noindent
\textbf{Acknowledgment} \ \ It is our pleasure to thank Takaaki Kagawa for his reading the first version of our paper and making various comments including the references to his pioneer works. We would like to thank Masanari Kida and Yuichiro Taguchi who gave us some useful advice. Our thanks also go to Iwao Kimura and Denis Simon who pointed out that there were some bugs of programs what we used.

\section{Setup}

In this section, we introduce the strategy to prove our results. Henceforth, we assume that the class number of $K_m$ is 1 and every elliptic curve $E$ with everywhere good reduction over $K_m$ has no $K_m$-rational point of order $2$ because Comalada \cite{Com90} determines all admissible curves ($=$ elliptic curves having good reduction everywhere and a $K_m$-rational point of order $2$) defined over $K_m$ with $m \le 100$ and such curves do not exist over $K_m$ which we consider in this paper. First we use the following result:

\begin{prop}[Setzer \cite{Set78}]
\textit{Let $E$ be an elliptic curve over $K_m$. If the class number of $K_m$ is prime to $6$ then $E$ has a global minimal model.}
\end{prop}


Let $E$ be an elliptic curve with everywhere good reduction over $K_m$. By Proposition 2.1, $E$ has a global minimal model
\[
E: \ y^2 + a_1 xy + a_3y = x^3 + a_2 x^2 + a_4x +a_6
\]
with coefficients $a_i \in \mathcal{O}_{K_m}$ $(i=1,2,3,4,6)$. The discriminant of $E$ (denoted by $\Delta(E)$) is
\[
\Delta(E) = \frac{c_4^3-c_6^2}{1728}
\]
where $c_4,c_6 \in \mathcal{O}_{K_m}$ are, as in \cite{Sil09} (Chapter III, p.42), written as polynomials in the $a_i$'s with integer coefficients. 
Moreover, the following conditions are equivalent (cf.\ \cite{Sil09}, Chapter VII, Prop.\ 5.1):
\begin{itemize}
\item $E$ has everywhere good reduction over $K_m$,
\item $\Delta(E) \in \mathcal{O}_{K_m}^{\times}$.
\end{itemize}
In our case, all elements of $\mathcal{O}_{K_m}^{\times}$ are written in the form $\pm \varepsilon^n$ where $\varepsilon$ is a fundamental unit of $K_m$ (let us fix $\varepsilon$ for each $m$).
Thus to determine the elliptic curves with everywhere good reduction over $K_m$, we shall compute the sets
\[
E_n^{\pm} (\mathcal{O}_{K_m}) = \left\{ (x,y) \in \mathcal{O}_{K_m} \times \mathcal{O}_{K_m} \mid \ y^2 = x^3 \pm 1728 \varepsilon^n  \right\}, \ \ 0 \le n < 12.
\]
However, the set of coefficients $(a_1, a_2, a_3, a_4, a_6) \in \mathcal{O}_{K_m}^{\oplus 5}$, which gives rise to $(c_4,c_6) \in \mathcal{O}_{K_m}^{\oplus 2}$, does not necessarily exist. Therefore, we check whether the curve
\begin{eqnarray}
E_C: \ y^2 = x^3 - 27 c_4 x -54 c_6 ,
\end{eqnarray}
which is isomorphic to $E$ over $K_m$, has trivial conductor for each $(c_4,c_6) \in E_n^{\pm} (\mathcal{O}_{K_m})$.



Actually, it is very hard to compute all $E_n^{\pm} (\mathcal{O}_{K_m})$ because of the limitation of efficiency of equipments. To reduce the amount of computation, we show that some values of $n$ are irrelevant by using Kagawa's results. In \cite{Kag11}, Kagawa shows a criterion whether the discriminant of an elliptic curve with everywhere good reduction over $K_m$ is a cube in $K_m$:

\begin{lem}[\cite{Kag11}, Prop.\ 1]
\textit{If the following five conditions hold, then the discriminant of every elliptic curve with everywhere good reduction over $K_m$ is a cube in $K_m$:}
\begin{enumerate}
\item \textit{The class number of $K_m$ is prime to $6$;}
\item \textit{$K_m/\mathbb{Q}$ is unramified at $3$;}
\item \textit{The class number of $K_m (\sqrt{-3})$ is prime to $3$;}
\item \textit{The class number of $K_m (\sqrt[3]{\varepsilon})$ is odd;}
\item \textit{For some prime ideal $\mathfrak{p}$ of $K_m$ dividing $3$, the congruence $X^3 \equiv \varepsilon$ $(\mathrm{mod} \ \mathfrak{p}^2)$ does not have a solution in $\mathcal{O}_{K_m}$.}
\end{enumerate}
\end{lem}

Using the criterion, Kagawa shows the following:

\begin{lem}[\cite{KagTC}]
\textit{If $m=46$ or $59$, every elliptic curve with everywhere good reduction over $K_m$ has a global minimal model whose discriminant is a cube in $K_m$.}
\end{lem}

\noindent
Therefore, we have $\Delta(E) =  \pm \varepsilon^{3n}$ for some $n \in \mathbb{Z}$.

By applying the next lemma, we can further discard some cases:

\begin{lem}[\cite{Kag11}, Prop.\ 4]
\textit{Let $E$ be an elliptic curve defined over $K_m$. If $E$ has good reduction outside $2$ and has no $K_m$-rational point of order $2$, then $K_m (E[2]) / K_m (\sqrt{\Delta (E)})$ is a cyclic cubic extension unramified outside $2$. In particular, the ray class number of $K_m (\sqrt{\Delta (E)})$ modulo $\prod_{\mathfrak{p}|2} \mathfrak{p}$ is a multiple of $3$.}
\end{lem}

\noindent
Note that $K_m (\sqrt{\Delta (E)})$ is either $K_m$, $K_m (\sqrt{-1})$ or $K_m (\sqrt{\pm \varepsilon})$.
Thus we compute the ray class number of $K_m$ modulo $\prod_{\mathfrak{p}|2} \mathfrak{p}$. The following computations are carried out by using Pari/GP \cite{Pari} (Same type results were obtained in \cite{KagPh} by using KASH \cite{KASH}). The bold-faced numbers in this table are the ones divisible by $3$.

\begin{center}
\begin{tabular}{|c||c|c|c|c|}
\hline
$m$ & $K_m$ & $K_m (\sqrt{-1})$ & $K_m (\sqrt{\varepsilon})$ & $K_m (\sqrt{- \varepsilon})$  \\ \hline \hline
$43$ & 1 & \textbf{3} & 10 & 1  \\ \hline
$46$ & 1 & 4 & 1 & \textbf{3}  \\ \hline
$59$ & 1 & \textbf{9} & \textbf{6} & 1 \\ \hline
\end{tabular}

\vspace{2mm}

Table 1. Ray class number of $K_m (\sqrt{\Delta(E)})$ $(m=43,46,59)$ modulo $\prod_{\mathfrak{p}|2} \mathfrak{p}$ 
\end{center}

\noindent
As a result, if $m=46$ the discriminant $\Delta(E)$ is $- \varepsilon^{6n+3}$ and if $m=59$ the discriminant is $- \varepsilon^{6n}$ or $\varepsilon^{6n+3}$.
We can conclude that it is enough to determine $E_3^+ (\mathcal{O}_{K_{46}})$, $E_0^+ (\mathcal{O}_{K_{59}})$ and $E_3^- (\mathcal{O}_{K_{59}})$ to prove Theorem 1.2 for $m=46$ and $59$.

However, the case $m=43$ remains because some of the conditions in Lemma 2.2 do not hold. In this case, it is known that the discriminant is $- \varepsilon^{2n}$ (cf.\ \cite{KagPh} and Lem.\ 2.3) so we need to compute the three sets, $E_0^+ (\mathcal{O}_{K_{43}})$, $E_2^+ (\mathcal{O}_{K_{43}})$ and $E_4^+ (\mathcal{O}_{K_{43}})$.

\section{Results of the computation}

\subsection{Computing Mordell-Weil groups and integral points}

To compute $E_n^{\pm} (\mathcal{O}_{K_m})$, we first compute the Mordell-Weil group $E_n^{\pm} (K_m)$. It is decomposed into a direct-sum of $E_n^{\pm} (K_m)_{\mathrm{tors}}$ (torsion part) and $E_n^{\pm} (K_m)_{\mathrm{free}}$ (free part, which is not canonical). The torsion part can be determined by observing reduction at good primes and decomposition of division polynomials. On the other hand, the free part can be computed by applying two-descent and infinite descent (the process of decompression from $E_n^{\pm} (K_m) / 2 E_n^{\pm} (K_m)$ to $E_n^{\pm} (K_m)$).

\begin{prop}
\textit{A basis of $E^{\pm}_n (K_m)$ is as follows:}
\begin{enumerate}
\item (Case $m=43$)
\begin{enumerate}
\item \textit{$E_0^+ (K_{43}) \simeq \mathbb{Z} \oplus \mathbb{Z} / 2 \mathbb{Z}$ and a basis is $\left\{ T_{43}, P_{43A} \right\}$ where $T_{43}=(-12,0)$ is $2$-torsion and}
\[
P_{43A}=\left( -\frac{104}{9} , - \frac{56}{27} \sqrt{43} \right) 
\]
\textit{is a generator of the free-part.}
\item \textit{$E_2^+ (K_{43}) \simeq \mathbb{Z}$ and a basis is $\left\{ P_{43B} \right\}$ where}
\[
P_{43B}=\left( 3200-488 \sqrt{43} , 294088-44848 \sqrt{43} \right) .
\]
\item \textit{$E_4^+ (K_{43}) \simeq \mathbb{Z}$ and a basis is $\left\{ P_{43C} \right\}$ where}
\[
P_{43C}=\left( -727456 + 110936 \sqrt{43} , 496115392 - 75656888 \sqrt{43} \right) .
\]
\end{enumerate}

\item (Case $m=46$) \textit{$E_3^+ (K_{46})$ is isomorphic to $\mathbb{Z} \oplus \mathbb{Z} / 2 \mathbb{Z}$ with a basis $\left\{ T_{46}, P_{46} \right\}$ where $T_{46}=(-12 \varepsilon,0)$ $(\varepsilon = 24335 + 3588 \sqrt{46})$ is $2$-torsion and}
\[
P_{46}=\left( \frac{1044823225}{6084} + \frac{987505}{39} \sqrt{46} , \frac{116177050458217}{474552} + \frac{73202442649}{2028} \sqrt{46} \right) 
\]
\textit{is a generator of the free-part.}
\item (Case $m=59$)
\begin{enumerate}
\item \textit{$E_0^+ (K_{59}) \simeq \mathbb{Z} \oplus \mathbb{Z} / 2 \mathbb{Z}$ and a basis is $\left\{ T_{59A}, P_{59A} \right\}$ where $T_{59A}=(-12,0)$ is $2$-torsion and}
\[
P_{59A}=\left( - \frac{133}{16} , \frac{283}{64} \sqrt{59} \right) 
\]
\textit{is a generator of the free-part.}
\item \textit{$E_3^- (K_{59}) \simeq \mathbb{Z}^{\oplus 2} \oplus \mathbb{Z} / 2 \mathbb{Z}$ and a basis is $\left\{ T_{59B}, P_{59B}, P_{59C} \right\}$ where $T_{59B}=(12 \varepsilon ,0)$ $(\varepsilon = -530 + 69 \sqrt{59})$ is $2$-torsion and}
\[
P_{59B}=\left( 9275 - \frac{2415}{2} \sqrt{59} , -\frac{5810733}{4} + \frac{756493}{4} \sqrt{59} \right)
\]
\textit{and}
\[
P_{59C}=\left( \frac{50000200}{59} - \frac{6509460}{59} \sqrt{59} , \frac{65094772968}{59} - \frac{500002437752}{3481} \sqrt{59} \right) 
\]
\textit{are generators of the free-part.}
\end{enumerate}
\end{enumerate}
\end{prop}

Here we used Denis Simon's two-descent program (cf.\ \cite{Sim02}) on Pari-GP \cite{Pari}. To compute some related data efficiently, we executed the Pari-GP program on Sage \cite{Sage} as a built-in software.

\vspace{3mm}

\noindent
\textbf{Warning:} Simon's two-descent program is also available as a Sage's built-in function that does not require the Pari-GP platform. However, this function has fatal bugs (errors) that come from the same bugs in the previous edition of Simon's original (Pari-GP) program. For the Pari-GP platform, this problem has already been fixed by himself completely, but it is not yet for the Sage platform.

\vspace{3mm}

To compute the subset $E_n^{\pm} (\mathcal{O}_{K_m})$ of integral points in $E_n^{\pm} (K_m)$, we use the method of elliptic logarithm to compute the linear form:
\[
P = \sum_{i=1}^r m_i P_i + n T \in E_n^{\pm} (\mathcal{O}_{K_m}) \ \ \ \ \ (m_1, ... , m_r, n \in \mathbb{Z})
\]
where $P_i$'s and $T$ are generators of the free part and the torsion part. Moreover, the maximum of the absolute values of the coefficients of the linear form
\[
M := \max \left\{ |m_1|, ..., |m_r| , |n| \right\}
\]
can be bounded using the LLL-algorithm (by Lenstra-Lenstra-Lov\'{a}sz, cf.\ \cite{Smart}).

\begin{prop}
\textit{The set of integral points $E_n^{\pm} (\mathcal{O}_{K_m})$ is as follows:}
\begin{enumerate}
\item (Case $m=43$)
\begin{enumerate}
\item $E_0^+ (\mathcal{O}_{K_{43}}) = \left\{ O, T_{43}, T_{43} \pm P_{43A} \right\}$,
\item $E_2^+ (\mathcal{O}_{K_{43}}) = \left\{ O, \pm P_{43B}, \pm 2P_{43B}  \right\}$,
\item $E_4^+ (\mathcal{O}_{K_{43}}) = \left\{ O, \pm P_{43C}, \pm 2P_{43C}  \right\}$.
\end{enumerate}
\item (Case $m=46$) $E_3^+ (\mathcal{O}_{K_{46}}) = \left\{ O, T_{46} \right\}$.
\item (Case $m=59$)
\begin{enumerate}
\item $E_0^+ (\mathcal{O}_{K_{59}}) = \left\{ O, T_{59A} \right\}$,
\item $E_3^- (\mathcal{O}_{K_{59}}) = \left\{ O, T_{59B} \right\}$.
\end{enumerate}
\end{enumerate}
\end{prop}

Finally, we compute that the elliptic curve $(1)$ has trivial conductor. As a result, there are no pairs $(c_4,c_6) \in E_n^{\pm} (\mathcal{O}_{K_m})$ for which $(1)$ has trivial conductor. Therefore, the non-existence of the curves follows.

\vspace{4mm}

In the same way, we can prove Theorem 1.3. For $m=62,67$ and $71$, by the assumption of cubic discriminants, it is enough to determine $E_{3n}^{\pm} (\mathcal{O}_{K_m})$ ($n \in \mathbb{Z}$). To apply Lemma 2.4, we compute the ray class number of $K_m$ modulo $\prod_{\mathfrak{p}|2} \mathfrak{p}$.

\begin{center}
\begin{tabular}{|c||c|c|c|c|}
\hline
$m$ & $K_m$ & $K_m (\sqrt{-1})$ & $K_m (\sqrt{\varepsilon})$ & $K_m (\sqrt{- \varepsilon})$  \\ \hline \hline
$62$ & 1 & 8 & \textbf{3} & 1  \\ \hline
$67$ & 1 & \textbf{3} & 14 & 1 \\ \hline
$71$ & 1 & 7 & \textbf{3} & 4 \\ \hline
\end{tabular}

\vspace{2mm}

Table 2. Ray class number of $K_m (\sqrt{\Delta(E)})$ $(m=62,67,71)$ modulo $\prod_{\mathfrak{p}|2} \mathfrak{p}$ 
\end{center}

\noindent
Finally it is enough to determine $E_3^- (\mathcal{O}_{K_{62}})$, $E_0^+ (\mathcal{O}_{K_{67}})$ and $E_3^- (\mathcal{O}_{K_{71}})$. Here is the result of computing Mordell-Weil groups and sets of integral points.

\begin{prop}
\textit{A basis of $E^{\pm}_n (K_m)$ and the set of integral points $E^{\pm}_n (\mathcal{O}_{K_m})$ are as follows:}
\begin{enumerate}
\item (Case $m=62$) \textit{$E_3^- (K_{62})$ is isomorphic to $\mathbb{Z} \oplus \mathbb{Z} / 2 \mathbb{Z}$ with a basis $\left\{ T_{62}, P_{62} \right\}$ where $T_{62}=(12 \varepsilon ,0)$ $(\varepsilon = -63 + 8 \sqrt{62})$ is $2$-torsion and}
\[
P_{62}=\left( \frac{30492}{25} - \frac{3872}{25} \sqrt{62} , -\frac{8377936}{125} + 8512 \sqrt{62} \right)
\]
\textit{is a generator of the free-part. The set of integral points is}
\[
E_3^- (\mathcal{O}_{K_{62}}) = \left\{ O,T_{62}  \right\} .
\]

\item (Case $m=67$) \textit{$E_0^+ (K_{67})$ is isomorphic to $\mathbb{Z} \oplus \mathbb{Z} / 2 \mathbb{Z}$ with a basis $\left\{ T_{67}, P_{67} \right\}$ where $T_{67}=(-12,0)$ is $2$-torsion and}
\[
P_{67}=\left( -\frac{584}{49} , \frac{248}{343} \sqrt{67} \right) 
\]
\textit{is a generator of the free-part. The set of integral points is}
\[
E_0^+ (\mathcal{O}_{K_{67}}) = \left\{ O,T_{67}, T_{67} \pm P_{67}   \right\} .
\]
\item (Case $m=71$) \textit{$E_3^- (K_{71})$ is isomorphic to $\mathbb{Z}^{\oplus 2} \oplus \mathbb{Z} / 2 \mathbb{Z}$ with a basis $\left\{ T_{71}, P_{71A}, P_{71B} \right\}$ where $T_{71}=(12 \varepsilon ,0)$ $(\varepsilon = 3480 + 413 \sqrt{71})$ is $2$-torsion and}
\[
P_{71A}= \left( 165300+ \frac{39235}{2} \sqrt{71} , \frac{377098253}{4} + \frac{44753329}{4} \sqrt{71}  \right)
\]
\textit{and}
\[
P_{71B}= \left( \frac{1560462848}{3025} + \frac{185192868}{3025} \sqrt{71}, - \frac{87152513410872}{166375} - \frac{10343100438152}{166375} \sqrt{71} \right)
\]
\textit{are generators of the free-part. The set of integral points is}
\begin{center}
$E_3^- (\mathcal{O}_{K_{71}}) = \left\{ O, T_{71}, \pm P_{71A} \mp P_{71B} , T_{71} \pm P_{71A} \mp P_{71B}  \right\}$. \\
(\textit{double sign in the same order})
\end{center}
\end{enumerate}
\textit{Moreover, there are no pairs $(c_4,c_6) \in E_n^{\pm} (\mathcal{O}_{K_m})$ for which the elliptic curve $(1)$ has trivial conductor.}
\end{prop}

\subsection{Trial of computation}

In this subsection, we show examples of computation times of running Simon's two-descent program to compute Mordell-Weil groups. Simon's two-descent is mainly controlled by four parameters:
\begin{itemize}
\item \texttt{lim1}: limit on trivial points on binary quartic forms (``quartics'' for short),
\item \texttt{lim3}: limit on points on ELS (everywhere locally solvable) quartics,
\item \texttt{limtriv}: limit on trivial points on elliptic curve,
\item \texttt{limbigprime}: distinguish between small and large prime numbers to use probabilistic tests for large primes,
\end{itemize}
and there are some supplemental parameters (\texttt{maxprob},\texttt{bigint},\texttt{nbideaux}, etc.). 

Now we fix the set of main parameters (\texttt{lim1},\texttt{lim3},\texttt{limtriv},\texttt{limbigprime})$=(40,60,40,30)$ that were chosen to compute the case $m=43$. The total running times of these computations are as follows:

\begin{center}
\begin{tabular}{|c||c|c|c|c|c|}
\hline
$m$ & $E_n^{\pm}$ & desired & actual & CPU time (sec.) & S/F \\ \hline \hline
\multirow{3}{*}{$43$} & $E_0^+$ & 1 & 1 & 570.168 & success  \\
\cline{2-6} & $E_2^+$ & 1 & 1 & 120.916 & success  \\
\cline{2-6} & $E_4^+$ & 1 & 1 & 112.554 & success \\ \hline
$46$ & $E_3^+$ & 1 & 1 & 670.117 & success  \\ \hline
\multirow{2}{*}{$59$} & $E_0^+$ & 1 & 1 & 195.500 & success  \\
\cline{2-6} & $E_3^-$ & 2 & 1 & 300.582 & failure \\ \hline
$62$ & $E_3^-$ & 1 & 1 & 317.216 & success  \\ \hline
$67$ & $E_0^+$ & 1 & 0 & 976.785 & failure  \\ \hline
$71$ & $E_3^-$ & 2 & 2 & 279.413 & success \\ \hline
\end{tabular}

\vspace{1mm}

Table 3. Rank of $E_n^{\pm} (K_m)$ with computation time \\ (Intel Core$^{\mathrm{TM}}$ i5 processor (3.30GHz, dual core) and 4.0GB memory)
\end{center}

As above, our trials failed for two cases due to the difficulty in searching for points on these curves $E_n^{\pm}$. Thus we need to change these parameters to get our results. As a result, we succeed in computing the case $E_3^- (\mathcal{O}_{K_{59}})$ and $E_0^+ (\mathcal{O}_{K_{67}})$ with the set of parameters (\texttt{lim1},\texttt{lim3},\texttt{limtriv},\texttt{limbigprime})$=(70,80,150,80)$ but we need a lot of time for the computation.

\begin{center}
\begin{tabular}{|c||c|c|c|c|c|}
\hline
$m$ & $E_n^{\pm}$ & desired & actual & CPU time (sec.) & S/F \\ \hline \hline
$59$ & $E_3^-$ & 2 & 2 & 2192.911 & success  \\ \hline
$67$ & $E_0^+$ & 1 & 1 & 6186.093 & success  \\ \hline
\end{tabular}

\vspace{1mm}

Table 4. Rank of $E_n^{\pm} (K_m)$ with computation time \\ for the case $E_3^- (\mathcal{O}_{K_{59}})$ and $E_0^+ (\mathcal{O}_{K_{67}})$  \\ (Intel Core$^{\mathrm{TM}}$ i5 processor (3.30GHz, dual core) and 4.0GB memory)
\end{center}


\vspace{5mm}

\noindent
Shun'ichi Yokoyama and Yu Shimasaki\\
Graduate School of Mathematics, Kyushu University\\
744 Motooka, Nishi-ku, Fukuoka, 819-0395, Japan\\
E-mail address:\\
\texttt{s-yokoyama@math.kyushu-u.ac.jp}\\
\texttt{y-shimasaki@math.kyushu-u.ac.jp}

\end{document}